\documentclass[12pt,twoside]{amsart}
\usepackage{amsmath, amsthm, amscd, amsfonts, amssymb, graphicx}
\usepackage{enumerate}
\usepackage[colorlinks=true,
linkcolor=blue,
urlcolor=cyan,
citecolor=red]{hyperref}
\usepackage{mathrsfs}
\newtheorem{theorem}{Theorem}[section]
\newtheorem{lemma}{Lemma}[section]
\newtheorem{remark}{Remark}[section]

\newtheorem{corollary}{Corollary}[section]

\numberwithin{equation}{section}

\begin{document}

\title{Refined Uncertainty Relation for q-Commutator}
\author{Kenjiro Yanagi}

\begin{abstract}
We show that Robertson uncertainty relation can be refined for q-commutator which is defined by $[A,B]_q = AB-qBA$, where 
$A, B$ are self-adjoint operators and $q  \in \mathbb{R}$. The coefficient is represented by the eigenvalues of state $\rho$. 
\end{abstract}
\maketitle
\pagestyle{myheadings}
\markboth{\centerline {}}
{\centerline {}}
\bigskip
\bigskip
\section{Introduction}

It is well known that Robertson uncertainty relation is represented by the commutator in the following. 

\begin{theorem}
Let $\rho$ be a density operator and $A, B$ be self-adjoint operators.  
Then the following uncertainty relation holds.
$$
\frac{1}{4} |Tr[\rho [A,B]]|^2 \leq V_{\rho}(A) V_{\rho}(B),  
$$
where $V_{\rho}(A) = Tr[\rho A^2]-(Tr[\rho A])^2 = Tr[\rho A_0^2]$ and $A_0 = A-Tr[\rho A]I$ .
\label{th:theorem1.1}
\end{theorem}

\noindent
It is easy to show the following uncertainty relation for q-commutator. 

\begin{theorem}
Let $\rho$ be a density operator and $A, B$ be self-adjoint operators.  Then the following uncertainty relation holds for $q \in \mathbb{R}$.
$$
\frac{1}{(1+|q|)^2} |Tr[\rho [A_0,B_0]_{|q|}]|^2 \leq V_{\rho}(A) V_{\rho}(B). 
$$
\label{th:theorem1.2}
\end{theorem}

\noindent
In this paper we have an  extension of Theorem \ref{th:theorem1.2} by using the maximal eigenvalue and minimal eigenvaue 
of $\rho$ under finite dimensional case.


\section{q-commutator}

Let $M_n(\mathbb{C}), M_{n,sa}(\mathbb{C})$ and $M_{n,1}(\mathbb{C})$ be the set of all $n$ dimensional complex matrices, 
the set of all $n$ dimensional self-adjoint matrices and the set of all $n$ dimensional density matrices. For $A, B \in M_{n,sa}(\mathbb{C})$, 
we define commutator and anti-commutator by $[A,B] = AB-BA$ and $\{ A,B \} = AB+BA$, respectively. 
Q-commutator and q-anti-commutator are defined by $[A,B]_q = AB-qBA$ and 
$\{ A,B \}_q = AB+qBA$, respectively, where $q \in \mathbb{R}$ and $A, B \in M_{n,sa}(\mathbb{C})$. 
Q-commutator is a generalization of commutator $[A,B]$. Now we state the main theorem. 

\begin{theorem}
Let $\rho = \sum_{i=1}^n \lambda_i |\phi_i \rangle \langle \phi_i|$ be the spectral decomposition satisfying 
$0 \leq \lambda_1 \leq \lambda_2 \leq \cdots \leq \lambda_n$. Then we have the following refined uncertainty relations. 
\begin{enumerate}
\item[(1)]  If $|q| \leq 1$, then 
$$
\frac{(\lambda_n+|q|\lambda_1)^2}{(1+|q|)^2(\lambda_n-|q|\lambda_1)^2} |Tr[\rho [A_0,B_0]_{|q|}]|^2 \leq V_{\rho}(A) V_{\rho}(B).
$$
\item[(2)] If $|q| > 1$ then 
\begin{eqnarray*}
&   & \frac{q^2(|q|\lambda_n+\lambda_1)^2}{(1+|q|)^2(|q|\lambda_n-\lambda_1)^2} |Tr[\rho [A_0,B_0]_{1/|q|}]|^2 \\
& = & \frac{(|q|\lambda_n+\lambda_1)^2}{(1+|q|)^2(|q|\lambda_n-\lambda_1)^2} |Tr[\rho [B_0,A_0]_{|q|}]|^2 \leq V_{\rho}(A) V_{\rho}(B).
\end{eqnarray*}
\end{enumerate}
\label{th:theorem2.1}
\end{theorem}

\noindent
We need some lemmas in order to prove the main result.

\begin{lemma} 
If $|q| \leq 1$, then $G(t) = \left( \frac{t-|q|}{t+|q|} \right)^2$ is monotone increasing for $t \geq 1$. 
\label{lem:lemma2.1}
\end{lemma}

\noindent
{\bf Proof.}  $G^{'}(t) = \frac{4|q|(t-|q|)}{(t+|q|)^3} \geq 0$. \ \hfill $\Box$

\begin{lemma}
If $|q| \leq 1$, then $F(t) = (1+|q|t)^2(t-|q|)^2-(1-|q|t)^2(t+|q|)^2 \geq 0$ for $t \geq 1$.  \\
That is 
$$
\left( \frac{1-|q|t}{1+|q|t} \right)^2 \leq \left( \frac{t-|q|}{t+|q|} \right)^2.
$$
\label{lem:lemma2.2}
\end{lemma}

\noindent
{\bf Proof.} $F^{'}(t) = 4|q|(1-|q|^2)(3t^2-1) \geq 0$. Since $F(1) = 0$, we have $F(t) \geq 0$ for $t \geq 1$ \ \hfill $\Box$

\vspace{0.5cm}
\noindent
{\bf Proof of Theorem \ref{th:theorem2.1}.}  Let $a_{ij} = \langle \phi_i|A_0|\phi_j \rangle$ and $b_{ji} = \langle \phi_j|B_0|\phi_i \rangle$. \\
$(1)$~Since $Tr[\rho [A_0,B_0]_{|q|}] = \sum_{i,j} (\lambda_i-|q|\lambda_j)a_{ij}b_{ji}$, 
$$
|Tr[\rho [A_0,B_0]_{|q|}]| \leq \sum_{i,j} |\lambda_i-|q|\lambda_j| |a_{ij}| |b_{ji}|.
$$
By Schwarz inequality 
\begin{eqnarray*}
|Tr[\rho [A_0,B_0]_{|q|}]|^2 & \leq & \left( \sum_{i,j} |\lambda_i-|q|\lambda_j| |a_{ij}| |b_{ji}| \right)^2 \\
& \leq & \sum_{i,j}|\lambda_i-|q|\lambda_j| |a_{ij}|^2 \sum_{i,j} |\lambda_i-|q|\lambda_j | |b_{ji}|^2.
\end{eqnarray*}
Then
\begin{eqnarray}
&   & \frac{1}{(\lambda_n-|q|\lambda_1)^2} |Tr[\rho [A_0,B_0]_{|q|}]|^2 \nonumber \\
& \leq & \sum_{i,j} \left| \frac{\lambda_i-|q|\lambda_j}{\lambda_n-|q|\lambda_1}\right| |a_{ij}|^2 \cdot 
\sum_{i,j} \left| \frac{\lambda_i-|q|\lambda_j}{\lambda_n-|q|\lambda_1}\right| |b_{ji}|^2 \nonumber \\
& = & \left\{  \sum_{i<j} \left| \frac{\lambda_i-|q|\lambda_j}{\lambda_n-|q|\lambda_1}\right| |a_{ij}|^2
+\sum_{i \geq j} \left|\frac{\lambda_i-|q|\lambda_j}{\lambda_n-|q|\lambda_1}\right| |a_{ij}|^2 \right\} \label{eq:num2.1} \\
&   & \times \left\{ \sum_{i<j} \left| \frac{\lambda_i-|q|\lambda_j}{\lambda_n-|q|\lambda_1}\right| |b_{ji}|^2
+\sum_{i \geq j} \left| \frac{\lambda_i-|q|\lambda_j}{\lambda_n-|q|\lambda_1}\right| |b_{ji}|^2 \right\} .  \nonumber
\end{eqnarray}
When $\lambda_i \geq \lambda_j$, it follows the following from Lemma \ref{lem:lemma2.1}
$$
G(\frac{\lambda_i}{\lambda_j}) = \left( \frac{\lambda_i / \lambda_j-|q|}{\lambda_i / \lambda_j+|q|} \right)^2 \leq \left( \frac{\lambda_n/\lambda_1-|q|}{\lambda_n/\lambda_1+|q|} \right)^2 = G(\frac{\lambda_n}{\lambda_1}).
$$
Because $\lambda_i/\lambda_j \leq \lambda_n/\lambda_1$. Then we have 
$$
\left( \frac{\lambda_i-|q|\lambda_j}{\lambda_n-|q|\lambda_1} \right)^2 \leq \left( \frac{\lambda_i+|q|\lambda_j}{\lambda_n+|q|\lambda_1} \right)^2.
$$
That is 
$$
\left| \frac{\lambda_i-|q|\lambda_j}{\lambda_n-|q|\lambda_1} \right| \leq \left| \frac{\lambda_i+|q|\lambda_j}{\lambda_n+|q|\lambda_1} \right|.
$$
When $\lambda_i < \lambda_j$, it follows the following from Lemma \ref{lem:lemma2.2}
$$
\left( \frac{1-|q|(\lambda_j/\lambda_i)}{1+|q|(\lambda_j/\lambda_i)} \right)^2 \leq \left( \frac{\lambda_j/\lambda_i-|q|}{\lambda_j/\lambda_i+|q|} \right)^2.
$$
By Lemma \ref{lem:lemma2.1} we have 
$$
G(\frac{\lambda_j}{\lambda_i}) = \left( \frac{\lambda_j / \lambda_i-|q|}{\lambda_ij/ \lambda_i+|q|} \right)^2 \leq \left( \frac{\lambda_n/\lambda_1-|q|}{\lambda_n/\lambda_1+|q|} \right)^2 = G(\frac{\lambda_n}{\lambda_1}).
$$
Because $\lambda_j/\lambda_i \leq \lambda_n/\lambda_1$. Then we have 
$$
\left( \frac{\lambda_i-|q|\lambda_j}{\lambda_i+|q|\lambda_j} \right)^2 \leq \left( \frac{\lambda_n-|q|\lambda_1}{\lambda_n+|q|\lambda_1} \right)^2.
$$
That is 
$$
\left| \frac{\lambda_i-|q|\lambda_j}{\lambda_n-|q|\lambda_1} \right| \leq \left| \frac{\lambda_i+|q|\lambda_j}{\lambda_n+|q|\lambda_1} \right|.
$$
By (\ref{eq:num2.1}) we have 
$$
\frac{1}{(\lambda_n-|q|\lambda_1)^2} |Tr[\rho [A_0,B_0]_{|q|}]|^2 \leq \sum_{i,j} \left| \frac{\lambda_i+|q|\lambda_j}{\lambda_n+|q|\lambda_1}\right| |a_{ij}|^2 \cdot 
\sum_{i,j} \left| \frac{\lambda_i+|q|\lambda_j}{\lambda_n+|q|\lambda_1}\right| |b_{ji}|^2.
$$
Since $|a_{ij}| = |a_{ji}|$, we have 
\begin{eqnarray*}
&   & \sum_{i,j} |\lambda_i+|q|\lambda_j| |a_{ij}|^2 = \sum_{i,j} \lambda_i |a_{ij}|^2 + |q|\sum_{i,j} \lambda_j |a_{ij}|^2 \\
& = & \sum_{i,j} \lambda_i |a_{ij}|^2 + |q|\sum_{i,j} \lambda_i |a_{ji}|^2 = \sum_{i,j} \lambda_i |a_{ij}|^2 + |q|\sum_{i,j} \lambda_i |a_{ij}|^2 \\
& = & (1+|q|)\sum_{i,j} \lambda_i |a_{ij}|^2 = (1+|q|)V_{\rho}(A).
\end{eqnarray*}
Then 
$$
\frac{(\lambda_n+|q|\lambda_1)^2}{(1+|q|)^2(\lambda_n-|q|\lambda_1)^2} |Tr[\rho [A_0,B_0]_{|q|}]|^2 \leq V_{\rho}(A) V_{\rho}(B).
$$

\vspace{0.5cm}
\noindent
$(2)$~Since $1/{|q|} < 1$, we have the following by (1).
\begin{eqnarray*}
&   & \frac{(\lambda_n+(1/{|q|})\lambda_1)^2}{(1+1/{|q|})^2(\lambda_n-(1/{|q|})\lambda_1)^2} |Tr[\rho [A_0,B_0]_{1/|q|}]|^2 \\
& = & \frac{q^2 (|q|\lambda_n+\lambda_1)^2}{(1+|q|)^2(|q|\lambda_n-\lambda_1)^2} |Tr[\rho [A_0,B_0]_{1/{|q|}}]|^2 \leq V_{\rho}(A) V_{\rho}(B).
\end{eqnarray*}
Since 
$$
Tr[\rho [A_0,B_0]_{1/{|q|}}] = -\frac{1}{|q|}Tr[\rho [B_0,A_0]_{|q|}]|,
$$
we have 
$$
\frac{(|q|\lambda_n+\lambda_1)^2}{(1+|q|)^2(|q|\lambda_n-\lambda_1)^2} |Tr[\rho [B_0,A_0]_{|q|}]|^2 \leq V_{\rho}(A) V_{\rho}(B).
$$
\ \hfill $\Box$

\vspace{0.5cm}
\noindent
We obtain the following corollary.

\begin{corollary}
We have the following refined uncertainty relations under five conditions.
\begin{enumerate}
\item[(1)]  If $0 < q \leq 1$, then 
$$
\frac{(\lambda_n+q\lambda_1)^2}{(1+q)^2(\lambda_n-q\lambda_1)^2} |Tr[\rho [A_0,B_0]_q]|^2 \leq V_{\rho}(A) V_{\rho}(B).
$$
\item[(2)] If $q > 1$ then 
\begin{eqnarray*}
&   & \frac{q^2(q \lambda_n+\lambda_1)^2}{(1+q)^2(q \lambda_n-\lambda_1)^2} |Tr[\rho [A_0,B_0]_{1/q}]|^2 \\
& = & \frac{(q \lambda_n+\lambda_1)^2}{(1+q)^2(q \lambda_n+\lambda_1)^2} |Tr[\rho[B_0,A_0]_q]|^2 \leq V_{\rho}(A) V_{\rho}(B).
\end{eqnarray*}
\item[(3)] If $q = 0$, then $|Tr[\rho A_0 B_0]|^2 \leq V_{\rho}(A)V_{\rho}(B)$. 
\item[(4)] If $-1 \leq q < 0$, then 
\begin{eqnarray*}
&   & \frac{(\lambda_n-q\lambda_1)^2}{(1-q)^2(\lambda_n+q\lambda_1)^2} |Tr[\rho [A_0,B_0]_{-q}]|^2 \\
& = & \frac{(\lambda_n-q\lambda_1)^2}{(1-q)^2(\lambda_n+q\lambda_1)^2} |Tr[\rho \{A_0,B_0 \}_q]|^2 \leq V_{\rho}(A) V_{\rho}(B).
\end{eqnarray*}
\item[(5)] If $q < -1$, then
\begin{eqnarray*}
&   & \frac{q^2(\lambda_1-q\lambda_n)^2}{(1-q)^2(\lambda_1+q\lambda_n)^2} |Tr[\rho [A_0,B_0]_{-1/q}]|^2 \\
& = & \frac{q^2(\lambda_1-q\lambda_n)^2}{(1-q)^2(\lambda_1+q\lambda_n)^2} |Tr[\rho \{A_0,B_0 \}_{1/q}]|^2 \\
& = & \frac{(\lambda_1-q\lambda_n)^2}{(1-q)^2(\lambda_1+q \lambda_n)^2} |Tr[\rho \{ B_0,A_0 \}_q]|^2 \leq V_{\rho}(A) V_{\rho}(B).
\end{eqnarray*}
\end{enumerate}
\label{cor:corollary2.1}
\end{corollary}

\vspace{0.5cm}
\noindent
We obtain the result by \cite{KiMaYa:re}. This is the case of $q = 1$ in Theorem \ref{th:theorem2.1} and Corollary \ref{cor:corollary2.1}.

\begin{corollary}
\begin{eqnarray*}
&   & \frac{(\lambda_n+\lambda_1)^2}{4(\lambda_n-\lambda_1)^2}|Tr[\rho [A_0,B_0]]|^2 \\
& = & \frac{(\lambda_n+\lambda_1)^2}{4(\lambda_n-\lambda_1)^2}|Tr[\rho [A,B]]|^2    \leq V_{\rho}(A) V_{\rho}(B).
\end{eqnarray*}
\label{cor:corollary2.2}
\end{corollary}

\begin{remark}
Let  $\rho = \frac{1}{n}I$ be the maximally mixed state. If $|q| < 1$, then 
$$
\frac{1}{(1-|q|)^2} |Tr[\rho[A_0,B_0]_{|q|}]|^2 \leq V_{\rho}(A) V_{\rho}(B).
$$
If $|q| > 1$, then 
$$
\frac{1}{(1-|q|)^2} |Tr[\rho[B_0,A_0]_{|q|}]|^2 \leq V_{\rho}(A) V_{\rho}(B).
$$
These are the same as $|Tr[A_0B_0]|^2 \leq Tr[A_0^2]Tr[B_0^2]$.
\label{re:remark2.1}
\end{remark}

\begin{remark}
Let $\rho$ be non-faithfull state. Since $\lambda_1 = 0$, we have the same result as Theorem \ref{th:theorem1.2}. 
Because in the case of $|q| \leq 1$, 
$$
\frac{1}{(1+|q|)^2} |Tr[\rho[A_0,B_0]_{|q|}]|^2 \leq V_{\rho}(A) V_{\rho}(B).
$$
Then we can rewrite the following
$$
\frac{1}{(1+|q|)^2} |Tr[\rho[B_0,A_0]_{|q|}]|^2 \leq V_{\rho}(A) V_{\rho}(B).
$$
Then in the case of $|q| > 1$, we have
$$
\frac{|q|}{(1+|q|)^2} |Tr[\rho[B_0,A_0]_{1/|q|}]|^2 \leq V_{\rho}(A) V_{\rho}(B).
$$
Therefore 
$$
\frac{1}{(1+|q|)^2} |Tr[\rho[A_0,B_0]_{|q|}]|^2 \leq V_{\rho}(A) V_{\rho}(B).
$$
\label{re:remark2.2}
\end{remark}

\vspace{0.5cm}
\section*{Acknowledgement}

The author would like to thank Prof. G.Kimura for giving me interresting results. 


{\tiny (K. Yanagi) Emeritus Professor of Yamaguchi University, 2-16-1, Tokiwadai, Ube, 755-8611, Japan}

{\tiny \textit{E-mail address:} yanagi@yamaguchi-u.ac.jp}

\end{document}